\documentclass[reqno,12pt,letterpaper]{amsart}
\usepackage{sdmacros, enumerate}

\newcommand{\RR}{{\mathbb R}}
\newcommand{\ZZ}{{\mathbb Z}}
\newcommand{\NN}{{\mathbb N}}
\newcommand{\CC}{{\mathbb C}}

\title[ON THE POLLICOTT\---RUELLE RESONANCES]%
{On the Pollicott\---Ruelle resonances}
\author{Joel Antonio\---V\'asquez}
\email{hello@joelantonio.me}
\address{}
\begin{document}
\begin{abstract}
  The purpose of this survey is to present the recent advances about the
  Pollicott\---Ruelle resonances.
\end{abstract}

\maketitle

%%%%%%%%%%%%%%%%%%%%%%%%%%%%%%%%%%%%%%%%%%%%%%%%%%%%%%%%%%%%%%%%%%%%%%%%%%%%%%%%
%%%%%%%%%%%%%%%%%%%%%%%%%%%%%%%%%%%%%%%%%%%%%%%%%%%%%%%%%%%%%%%%%%%%%%%%%%%%%%%%

\section{introduction}
\hspace{-1.4em}
Suppose that $M$ is a smooth compact manifold and let $\varphi: M \longrightarrow M$ be a 
smooth flow such that $M$ is considered a $\phi_{t}$\---invariant set for some $t > 0$.
The flow $\varphi_{t}$ is called Anosov flow if rather is a hyperbolic set for $\varphi_{t}$
in the sense of \cite[Definition 6.4.18]{KaHa} and is generated by some smooth
vector field $V$ such that $\varphi_{t} := e^{tV}$. Let $f, g \in C^{\infty}(M)$ be smooth
functions and $\mu$ a $\phi$\---invariant probability measure, then the correlation
function is defined as 
\begin{equation}
  \label{eq:correlation}
  \rho_{f, g}(t) = \int_{M} f(\varphi_{-t}(x))g(x) d\mu, \hspace{1em} \text{for any $x \in M$}.
\end{equation}

The power spectrum of (\ref{eq:correlation}) is the Fourier transform such that
\begin{equation}
  \label{eq:fourier}
  \widehat{\rho}_{f, g}(\lambda) = \int_{0}^{\infty} \rho_{f, g}(t)e^{i \lambda t} dt,
\end{equation}
is meromorphic for $|\Imag \lambda| > 0$. Thus, the asymtoptic behaviour of $\rho_{f, g}$ is controlled by the poles of the
extension $\widehat{\rho}_{f, g}$, such poles are known as the \textit{Pollicott\---Ruelle resonances}.
In other words, they are complex numbers, which describe
fine of decay of correlations for an Anosov Flow on a smooth compact manifold,
and were initially studied by M. Pollicott \cite{Po85, Po86} and D. Ruelle \cite{Ru86, Ru87}.  
From another point of view, the Pollicott\---Ruelle resonances are also the
singularities of the meromorphic extension of the Ruelle zeta function,
which was conjectured by S. Smale in 1967 \cite{Sm}. Such conjecture, has been proved
by Giulietti\---Liverani-Pollicot \cite{GiLiPo} for compact manifolds. Later, Arnoldi-Faure-Weich
\cite{AFW} defined resonances on open hyperbolic surfaces and Faure\---Tsujii \cite{FaTsb}
defined resonances for the Grassmanian bundle of an Anosov flow. Recently, Dyatlov\---Guillarmou
\cite{DyGu14} were able to define Pollicott\---Ruelle resonances for open hyperbolic systems
on a more general way compared to \cite{AFW, FaTsb} via a microlocal approach of 
Faure\---Sj\"ostrand \cite{FaSj} and Dyatlov\---Zworski \cite{DyZw13}, holding the results of 
\cite[\S 7]{Po86} and as a consequence, they were able to show that the Ruelle zeta function 
extends meromorphically to the entire complex plane.   \\ \\
\textbf{Structure of the survey.} Section \ts \ref{section-compact-case} focus on the
compact case, while Section \ts \ref{section-open-case}
focus on the open case.

\section{On the compact case}
\label{section-compact-case} 
\subsection{On the functional analysis proof}
In 2012, Giulietti\---Liverani\---Pollicott \cite{GiLiPo} showed the existence of Pollicott\---Ruelle
resonances for compact manifolds, such proof was given thought the Ruelle zeta function
for $C^{r}$ Anosov flows for $r > 2$ on a compact smooth orientable manifold, where they 
proved that for $C^{\infty}$ flows the zeta function is meromorphic on the entire
complex plane. %Firstly, let us review some previous work on the compact case.  
Based on the statement that from $C^{r}$ flows, we can obtain a strip in which $\zruelle(z)$ is meromorphic
of width unboundedly increasing with $r$ \cite{Fr}, such work was an expansion of
\cite{GoLi, BuLi, LiTs, Lia, Lib, BaLi}.
\subsubsection{Definitions}
In this subsection, we work with the following assumptions:
\begin{enumerate}[\bf(B1)]
  \label{prop-B}
  \ii $M$ is a $d$\---dimensional connected, compact and orientable $C^{\infty}$
  Riemannian manifold for some $d \in \ZZ_{+}$, $V$ is a $C^{\infty}$ nonvanishing vector
  field on $M$, and $\varphi_{t} = e^{tV}$ is the corresponding flow;
  \ii for each $x \in M$, there is a splitting
  \begin{equation}
    \label{eq:splitting}
    T_{x}M = E_{s}(x) \oplus E_{0}(x) \oplus E_{u}(x),
  \end{equation}
  where $E_{0}$ is the one\---dimensional subspace tangent to the flow, such that for
  some constants $C, \gamma > 0$
   \begin{equation}\label{eq:hyper} 
     \begin{array}{lll}
       & |d\varphi_{t}(x) \cdot v| \leq  Ce^{-\gamma |t|}|v| \quad&\text{ if } t \geq 0,
       v \in E_{s} ;  \\
       & |d\varphi_{t}(x) \cdot v| \leq  Ce^{-\gamma |t|}|v| \quad & \text{ if } t \leq 0, 
       v \in E_{u} ;  \\
       & \hspace{-0.7em} C^{-1}|v| \leq |d\varphi_{t}(v)| \leq C|v|&  \text{ if } t \in \RR, v \in E_{0}.
     \end{array}
   \end{equation}
  We denote $d_{s} = \dimM(E_{s})$ and $d_{u} = \dimM(E_{u})$ to distinguish the
  dimension of the stable and unstable subspaces, respectively.
\end{enumerate}  
In the context of the assumptions $\hyperref[prop-B]{(B1)}$\---$\hyperref[prop-B]{(B2)}$, we
define the Ruelle zeta function as related to the Riemann zeta function
(\ie $\zeta(z) = \prod_{p}(1 - p^{-z})^{-1}$), replacing $p$ by primitive closed orbits. Thus,
  \begin{equation}
    \label{eq:zruelle}
    \zruelle(z) = \prod_{\tau \in \Tau_{p}} (1 - e^{-z \lambda(\tau)})^{-1}, \hspace{1em} z \in \CC,
  \end{equation}
where $\Tau_{p}$ denotes the set of prime orbits and $\lambda(\tau)$ denotes the 
periodic of the closed orbit $\tau$.  
According to \cite[\ts 10]{PaPo}, for weak mixing Anosov flows, the $\zruelle(z)$
is analytic and nonzero for $\Real(z) \geq \htop(\varphi_{1})$ apart for a single
pole at $z = \htop(\varphi_{1})$.  
In order to understand the whole case, some relevant definitions will be given below,
however, the assertions of some concepts will be cited and won't be part of the main proofs.
 
As part of the main definitions, Ruelle \cite{Ru76} related the transfer
operator with the dynamical Fredholm determinant and defined $\Ddet_{\ell}$ as
the dynamical determinants, which are functions defined from weighted periodic
orbit data of a \textit{differentiable dynamical system}.

\begin{defi}
  The dynamical (also known as Fredholm\---Ruelle) determinants are defined as
  \label{defi:dynamical-det}
  \begin{equation}
    \label{eq:dynamical-det}
    \Ddet_{\ell}(z) = \exp\left( - \sum_{\tau \in \Tau} \frac{\tr(\wedge^{\ell}
    (\Dhyp \varphi_{-\lambda(\tau)}))e^{-z\lambda(\tau)}}{\mu(\tau)\epsilon(\tau)
  |\det(\COne - \Dhyp \varphi_{-\lambda(\tau)})|} \right),
  \end{equation}
  where $\epsilon(\tau)$ is 1 if the flow preserves the orientation of $E_{s}$ along
  $\tau$ and -1 otherwise. More precisely,
  $\epsilon(\tau) = \sign(\det(D_{\varphi_{-\lambda(\tau)}} |_{E_{s}}))$.
\end{defi}  \hspace{-1.4em}
The symbol $\Dhyp \varphi_{-t}$ in Definition \ref{defi:dynamical-det}, indicates
the derivative of
the map induced by the local transverse sections to the orbit (one at $x$, the other
at $\varphi_{t}(x)$) and can be represented as a $(d-1)\times(d-1)$\---dimensional matrix.
By $\wedge^{\ell}A$ we mean the matrix associated to the standard $\ell$\---th exterior
product of $A$ \--- see more about dynamical determinants in \cite{Ba16, BaTs08}.  
Given any $\phi$\---invariant probability measure $\mu$ on $M$ with $h_{\mu}(\phi)$ being
the measure theoric entropy of $\phi_{1}$. The topological entropy $\htop(\phi_{1})$
can be defined by
\begin{equation}
  \label{eq:htop}
  \htop(\phi_{1}) \equiv \sup \{h_{\mu}(\phi): \mu \text{ is a $\phi$\---invariant
  probability measure}\}.
\end{equation}
Thus, given $0 \leq \ell \leq d-1, \tau \in \Tau$, we let
$$
    \chi_{\ell}(\tau) = \frac{\tr(\wedge^{\ell}(\Dhyp \varphi_{-\lambda(\tau)}))}{
    \epsilon(\tau)|\det(\COne - \Dhyp \varphi_{-\lambda(\tau)})|},
    $$
  in order to write Equation (\ref{eq:dynamical-det}) in a shorter way as
  \begin{equation}
    \label{eq:short-det}
    \Ddet_{\ell}(z) = \exp\left( -\sum_{\tau \in \Tau} \frac{\chi_{\ell}(\tau)}{\mu(\tau)}
    e^{-z\lambda(\tau)} \right). 
  \end{equation}
Now let us define $C^{r}$ sections of $\wedge^{\ell}(T^{*}M)$ as the space $\Omega_{v}^{\ell}(M)$
of $\ell$\---forms on $M$ for all $v, \ell \in \NN$.
\begin{defi}
  \label{defi:spaces-l-forms}
Let $\Omega_{0, v}^{\ell} \subset \Omega_{v}^{\ell}(M)$
be the subspace of forms null in the flow direction, such that
$$
\Omega_{0, v}^{\ell}(M) = \{h \in \Omega_{r}^{\ell}(M): h(V, \dots) = 0\}.
$$
\end{defi}  
For a detailed construction and proof of Definition \ref{defi:spaces-l-forms} \--- see \cite[\ts 3]{GiLiPo}.  
Gou\"ezel\---Liverani \cite{GoLi}, defined Banach spaces adapted to Anosov systems and was adapted by \cite{GiLiPo}
in the following way.

\begin{defi}
  \label{defi:bbanach}
  For all $p \in \NN, q \in \RR_{+}, \ell \in \{0, \dots, d-1\}$ we define the spaces
  $\Bbanach^{p, q, \ell}$ to be the closures of $\Omega_{0, r}^{\ell}(M)$ with respect
  to the norm $\norm{\cdot}_{-,p,q,\ell}$ and the spaces $\Bbanach_{+}^{p, q, \ell}$
  to be the closures of $\Omega_{0, r}^{\ell}(M)$ with respect to the norm 
  $\norm{\cdot}_{+,p,q,\ell}$.
\end{defi} \hspace{-1.4em}
The construction of the $\Bbanach$ space in Definition \ref{defi:bbanach} is detailed in \cite[\ts 3.2]{GiLiPo}.
The sum in Equation (\ref{eq:dynamical-det}) is well\---defined, provided $\Real(z)$ is large enough
and follows a product analogous \cite[Equation (2.5)]{GiLiPo} such that
\begin{equation}
  \label{eq:zruella-det}
  \prod_{\ell = 0}^{d-1} \Ddet_{\ell}(z)^{(-1)^{\ell + d_{s} + 1}} = \zruelle(z).
\end{equation}
Now, to take care of the $t \leq t_{0}$, we introduce the dynamical norm
$\norm{\cdot}_{p,q,\ell}$ \--- see \cite[\ts 4]{GiLiPo}. For each $h \in \Omega_{\tau}^{\ell}(M)$, we set
\begin{equation}
  \label{eq:ddet-h}
  \norm{h}_{p,q,\ell} = \sup_{s \leq t_{0}} \norm{\Lop_{s}^{(\ell)}h}_{p,q,\ell},
\end{equation}
where $\Lop_{t}^{(\ell)}$ is a linear operator such that 
$\Lop_{t}^{(\ell)}: \Omega_{0, r-1}^{\ell}(M) \longrightarrow \Omega_{0, r-1}^{\ell}(M)$,
  for some $t \in \RR_{+}$. Furthermore,
  \begin{equation}
    \label{eq:lop}
    \Lop_{t}^{(\ell)}h := \varphi_{-t}^{*}h,
  \end{equation}
  for some $h \in \Omega_{0, r-1}^{\ell}(M)$ \--- see more \cite[\ts 4.1]{GiLiPo}.
Thus, we define 
$$\Cbanach^{p,q,\ell} = \overline{\Omega_{0,r}^{\ell}}^{\norm{\cdot}_{p,q,\ell}}
\subset \Bbanach^{p,q,\ell}.
$$ 
Let $\lambda_{i, \ell}$ be the eigenvalues of $X^{(\ell)}$. Then for each $z \in B_{}(\xi,  
\rho_{p,q, \ell})$, we let 
\begin{equation}
  \label{eq:ddet-analytic}
  \Cbanach(\xi - z, \xi) = \left( \prod_{\lambda_{i, \ell} \in B_{}(\xi, \rho_{p, q, \ell})}  
  \frac{z - \lambda_{i, \ell}}{\xi - \lambda_{i, \ell}} \right) \psi(\xi, z),
\end{equation} 
where $\psi(\xi, z)$ is analytic and nonzero for $z \in B_{}(\xi, \rho_{p,q, \ell})$. Thus,
the Equation (\ref{eq:ddet-analytic}) shows that the poles of $\zruelle$ are a 
subset of the eigenvalues of the $X^{(\ell)}$. 
\begin{defi}
  \label{defi:flat-trace}
  Given an operator $A \in L(\Bbanach^{p, q, \ell}, \Bbanach^{p, q, \ell})$, we define the
  flat trace as
  \begin{equation}
    \label{eq:flat-trace}
    \ftrace(A) = \lim_{\epsilon \rightarrow 0} \int_{M} \sum_{\alpha, \overline{i}}
  \iprod{\omega_{\alpha, \overline{i}}}{A(j_{\epsilon, \alpha, \overline{i}, x})
}_{x} \omega_{M}(x),
  \end{equation}
  where $\omega_{\alpha, i}$ is the dual of 1\---forms such that
  $\omega_{\alpha, i}(\hat{e}_{\alpha, j}) = \delta_{i, j}$ and 
  $j_{\epsilon, \alpha, \overline{i}, x}(y)$ is defined in \cite[Equation (5.2)]{GiLiPo},
  then provided the limit exists.
\end{defi} \hspace{-1.4em}
For a detailed proof of Definition \ref{defi:flat-trace} \--- see \cite[\ts 5]{GiLiPo}.

\subsubsection{On the proof}
We start  stablishing in what region $\zruelle(z)$ is meromorphic. %Thus, we have the following lemma.

\begin{lemm}
  \label{lemm:meromorphic-region}
  For any $C^{r}$ Anosov flow $\varphi_{t}$ with $r > 2$, then $\zruelle(z)$
  is meromorphic in the region
  \begin{equation}
      \Real(z) > \htop(\varphi_{1}) - \frac{\overline{\lambda}}{2} \floor{\frac{r-1}{2}}
  \end{equation}
  where $\overline{\lambda}$ is determined by the Anosov splitting.
\end{lemm} \hspace{-1.4em}
Lemma \ref{lemm:meromorphic-region} follows by the study of dynamical
determinants (Definition \ref{defi:dynamical-det}). In fact, to study in what
region $\zruelle$ is meromorphic, we must study in what region the dynamical
determinants are so. Moreover, for $\xi, z \in \CC$ we let
\begin{equation}
  \label{eq:two-det}
  \Cdet(\xi, z) = \exp\left( -\sum_{n=1}^{\infty} \frac{\xi^{n}}{n!}
  \sum_{\tau \in \Tau} \frac{\chi_{\ell}(\tau)}{\mu(\tau)} \lambda(\tau)^{n}
  e^{-z \lambda(\tau)} \right).
\end{equation}  
Using (\ref{eq:two-det}), for $\Real(z)$ sufficiently large and $|\xi - z|$
sufficiently small, we can write

$$
\Cdet(\xi - z, \xi) = \exp\left( -\sum_{n=1}^{\infty} \frac{(\xi - z)^{n}}{n!}
  \sum_{\tau \in \Tau} \frac{\chi_{\ell}(\tau)}{\mu(\tau)} \lambda(\tau)^{n}
  e^{-\xi \lambda(\tau)} \right) 
$$

$$
\hspace{2em} = \exp\left( -\sum_{\tau \in \Tau} \frac{\chi_{\ell}(\tau)}{\mu(\tau)} 
(e^{-z\lambda(\tau)} - e^{-\xi \lambda(\tau)})\right)
$$ 
\begin{equation}
\label{eq:ddet-reduce}
\hspace{-11.4em} = \frac{\Ddet_{\ell}(z)}{\Ddet_{\ell}(\xi)}.
\end{equation} 

\begin{theo}
  \label{theo:analytic-nonzero}
  $\zruelle(z)$ is analytic for $\Real(z) > \htop(\varphi_{1})$ and nonzero for
  $\Real(z) > \max\{\htop(\varphi_{1}) - \frac{\overline{\lambda}}{2} \floor{\frac{r-1}{2}}, 
  \htop(\varphi_{1}) - \overline{\lambda}\}$. Furthermore, if the flow is topologically
  mixing then $\zruelle(z)$ has no poles on the line $\{\htop(\varphi_{1}) + ib\}_{b \in \RR}$ 
  apart from the single simple pole at $z = \htop(\varphi_{1})$.
\end{theo} \hspace{-1.4em}
By Theorem \ref{theo:analytic-nonzero}, $\zruelle(z)$ is meromorphic in the entire complex
plane for smooth geodesic flows on any manifold that asserts the assumptions 
$\hyperref[prop-B]{(B1)}$\---$\hyperref[prop-B]{(B2)}$. Moreover, $\zruelle(z)$ has no zeroes
or poles on the line $\{\htop(\varphi_{1} + ib)\}_{b \in \RR}$, except at $z = \htop(\varphi_{1})$
where $\zruelle(z)^{-1}$ has a simple zero.  
From here, \cite{GiLiPo} specializes to contact Anosov flows. Let $\lambda_{+} \geq 0$ such that
$\norm{D\varphi_{-t}}_{\infty} \leq C_{0}e^{\lambda_{+}t}$ for all $t \geq 0$.

\begin{theo}
  \label{theo:cflow-pole}
  For a contact Anosov flow $\varphi_{1} \in C^{r}$ where $r > 2$, with 
  $\frac{\overline{\lambda}}{\lambda_{+}} > \frac{1}{3}$ there exists $\tau_{*} > 0$ such
  that the Ruelle zeta function is analytic in $\{z \in \CC : \Real(z) \geq 
  \htop(\varphi_{1}) - \tau_{*} \}$ apart from a simple pole at $z = \htop(\varphi_{1})$.
\end{theo}
\proof Equation (\ref{eq:ddet-ftrace}) and Equation (\ref{eq:ddet-reduce})
show that the poles of $\zruelle(z)$ are a subset of the eigenvalues of $X^{(\ell)}$.

\begin{lemm}
  \label{lemm:ddet-analytic}
  For any $C^{r}$ Anosov flow $\varphi_{t}$ with $r > 2, \xi \in \CC$ and $z \in \Ddet_{\ell}(\xi)$, 
  then it is analytic and nonzero in the region 
  $\Real(\xi) > \htop(\varphi_{1}) - \overline{\lambda}  |d_{s} - \ell|$.
\end{lemm}
\proof 
Let $\Omega_{0, v}^{\ell}(M) \subset \Omega_{v}^{\ell}(M)$ as Definition \ref{defi:spaces-l-forms},
let $\Bbanach^{p, q, \ell}$ such that $p \in \NN$ and $q \in \RR_{+}$ as Definition
\ref{defi:bbanach} and let $\Lop_{t}^{(\ell)}(h)$ as Equation (\ref{eq:ddet-h})
for some $h \in \Omega_{0, v}^{ \ell}(M)$.  
By restricting the transfer operator $\Lop_{t}^{(\ell)}$ to the space $\Omega_{0, r}^{\ell}(M)$ we
mimic the action of the standard transfer operators on sections transverse to the flow.
The operators (\ref{eq:lop}) generalize the action of the transfer operator $\Lop_{t}$ on
the spaces $\Bbanach^{p, q}$. Thus,
\begin{equation}
  \label{eq:dynorm1}
  \norm{\Lop_{t}^{\ell}h}_{p, q, \ell} \leq C_{p, q}e^{\sigma_{\ell}t} \norm{h}_{p, q, \ell},
\end{equation}
\begin{equation}
  \label{eq:dynorm2}
  \norm{h}_{p, q, \ell} = \sup_{s \leq t_{0}} \norm{\Lop_{s}^{(\ell)}h}_{p, q, \ell}.
\end{equation}
By Equations (\ref{eq:dynorm1}) and (\ref{eq:dynorm2}) imply that for some $t < t_{0}$, then
$$
\norm{\Lop_{t}^{(\ell)}h}_{p, q, \ell} \leq \max\{\norm{h}_{p, q, \ell}, C_{p, q}e^{|\sigma_{\ell}|t_{0}}
\norm{h}_{p, q, \ell}\} \leq C_{p, q} \norm{h}_{p, q, \ell},
$$
while for $t \geq t_{0}$ the required inequality holds trivially. The boundedness of $\Lop_{t}^{(\ell)}$
follows.  
The second inequality follows directly from the above, for small times and Equation (\ref{eq:dynorm1})
for larger times.  
On $\Cbanach^{p, q, \ell}$ the operators $\Lop_{t}^{\ell}$ form a strongly continuous
semigroup with generators $X^{(\ell)}$ by the above. We consider the resolvent $R^{(\ell)}(z) = 
(z\COne - X^{(\ell)})^{-1}$, then we have the following Lemma.
\begin{lemm}
  \label{lemm:quasi-comp-op}
  $R^{(\ell)}(z)$ is a quasi\---compact operator on $\Cbanach^{p, q, \ell}$.
\end{lemm}
\proof This follows by \cite[Lemma 3.8]{GiLiPo}. \qed   \\ \\
Although the operator $X^{(\ell)}$ is an unbounded closed
operator on $\Cbanach^{p, q, \ell}$, we can access to its spectrum thanks to Lemma
\ref{lemm:quasi-comp-op}. Now, let us make $\ell = d_{s}$ and let $\cvol_{s}$ be a volume.
A form on $E_{s}$ normalized so that $\norm{\cvol_{s}} = 1$ and it is globally continuous.
Let $\pi_{s}(x) = T_{x}M \longrightarrow E_{s}(x)$ be the projections on $E_{s}(x)$ 
along $E_{u}(x) \oplus E_{0}(x)$ such that
$$
\omega_{s}(v_{1}, \dots, v_{d_{s}}) = \cvol_{s}(\pi_{s}v_{1}, \dots, \pi_{s}v_{d_{s}})
$$
by construction $\omega_{s} \in \Omega_{0, \cvol}^{d_{s}}$. Note that $\varphi_{-t}^{*}\omega_{s}
= J_{s}\varphi_{-t}$ where $J_{s}\varphi_{-t}$ is the Jacobian restricted to the stable
manifold. Note that, $\omega_{s, \varepsilon} = \mathbb{M}_{\varepsilon}\omega_{s}$, for
$\varepsilon$ small enough, we have $\iprod{\omega_{s, \varepsilon}}{\omega_{s}} \geq \frac{1}{2}$.
Hence, 

\begin{equation}
  \label{eq:omg-jac}
  \int_{W_{\alpha}, G} \iprod{\omega_{s, \varepsilon}}{\Lop_{t}\omega_{s}} \geq
  \int_{W_{\alpha}, G} \frac{J_{s}\varphi_{-t}}{2} \geq \Csharp \int_{W_{\alpha, G}} J_{W}\varphi_{-t}
  \geq \Csharp \vol(\varphi_{-t}W_{\alpha, G}).
\end{equation}
Then the Equation (\ref{eq:omg-jac}) implies that the spectral radius of $R^{(d_{s})}(a)$ on
$\Cbanach^{q, \cvol, d_{s}}$ is exactly $(a - \sigma_{d_{s}})^{-1}$. Thus

\begin{equation}
  \label{eq:per-spec}
  \lim_{n \rightarrow \infty} \frac{1}{n} \sum_{k=0}^{n-1}(a-\sigma_{d_{s}})^{k} R^{(d_{s})}(a)^{k}
  = \left\{ \begin{array}{rcl}
      \prod & \mbox{if}(a-\sigma_{d_{s}})^{-1} \in \sigma_{\Cbanach^{p, q, \ell}}(R^{(d_{s})}(a)), 
      0 & \mbox{otherwise},
    \end{array}\right.
\end{equation}
where $\prod$ is the eigenprojector on the associated eigenspace and the convergence takes
places in the strong operator topology of $L(\Cbanach^{p, q, d_{s}}, \Cbanach^{p, q, d_{s}})$.
Thus, by Equation (\ref{eq:omg-jac}),
$$
\int_{W_{\alpha, G}} \iprod{W_{s, \varepsilon}}{\prod W_{s}} > 0,
$$
we have that $\prod \neq 0$ and $(a - \sigma_{d_{s}})^{-1}$ belongs to the spectrum.
This implies, that if $\ell = d_{s}$, then $\htop(\varphi_{1})$ is an eigenvalue of $X$ and
if the flow is topologically transitive $\htop(\varphi_{1})$ is a simple eigenvalue. Moreover,
if the flow is topologically mixing, then $\htop(\varphi_{1})$ is the only eigenvalue
on the line $\{\htop(\varphi_{1}) + ib\}_{b \in \RR}$. Thus, this proves Lemma \ref{lemm:ddet-analytic}.
\qed   \\ \\
In the same time, Theorem \ref{theo:cflow-pole} follows by Lemma \ref{lemm:ddet-analytic}.
\qed

\begin{theo}[{\cite[Theorem 2]{Po85}}]
  \label{theo:po85}
  Let $\varphi_{t}: \Lambda \longrightarrow \Lambda$ be a weak\---mixing Axiom A flow,
  then the Fourier transform $\fou_{f, g}(z)$  has a meromorphic extension to a strip
  $|\mathscr{J}(z)| \leq \varepsilon$, which is analytic on the real line. Furthermore,
  $\fou_{f,g}(t)$ tends to zero exponentially fast (for all H\"older continuous functions
  $f, g: \Lambda \longrightarrow \RR$) only if $\zeta(s, F)$ has an analytic extension
  to some strip $R(s) > P(F) - \varepsilon$, except for the simple pole at $s = P(F)$.
\end{theo}

\begin{theo}[{\cite[Paley-Wiener Theorem]{ReSi}}]
  \label{theo:resi}
  Let $\rho$ be in $\mathscr{S}'(\RR)$. Suposse that $\fou$ is a function with an analytic
  continuation to the set $\{\zeta \ | \ \Imag \zeta < a \}$ for some $a > 0$. Suppose
  further that for each $\eta \in \RR^{n}$ with $|\eta| < a$, $\fou(\cdot + i\eta)
  \in L^{1}(\RR^{n})$ and for any $b < a, \sup_{|\eta| < b} \norm{\fou(\cdot + i\eta)}_{1}
  < \infty$. Then $\rho$ is a bounded continuous function and for any $b < a$, there is
  a constant $C_{b}$ such that
  \begin{equation}
    \label{eq:paley-wiener}
  |\rho(x)| \leq C_{b}e^{-b|x|}.
  \end{equation}
\end{theo}

\begin{theo}
  \label{theo:ddet-nonzero}
  Let $\xi$ as in Lemma \ref{lemm:ddet-analytic}, then the function $\Cdet_{\ell}(\xi-z, \xi)$
  is analytic and nonzero for $z$ in the region
  \begin{equation}
  \label{eq:ddet-nonzero1}
  |\xi-z| < \Real(\xi) - \htop(\varphi_{1}) + |d_{s} - \ell|\overline{\lambda}
  \end{equation}
  and analytic in $z$, in the region
  \begin{equation}
  \label{eq:ddet-nonzero2}
  |\xi-z| < \Real(\xi) - \htop(\varphi_{1}) + |d_{s} - \ell|\overline{\lambda} +
  \frac{\overline{\lambda}}{2} \floor{\frac{r-1}{2}}
  \end{equation}
\end{theo}
\proof We can  write the spectral decomposition $R^{(\ell)}(z) = P^{(\ell)}(z) + U^{(\ell)}(z)$
where $P^{(\ell)}(z)$ is a finite rank operator and $U^{(\ell)}(z)$ has spectral radius 
arbitraly close to $\rho_{ess}(R^{(\ell)}(z))$.  
Let $\ftrace(R^{(\ell)}(z)^{n}) < \infty$ as in Definition \ref{defi:flat-trace}, then it can be written as
\begin{equation}
  \label{eq:ftrace2}
  \ftrace(R^{(\ell)}(z)^{n}) = \frac{1}{(n-1)!} \sum_{\tau \in \Tau}
  \frac{\chi_{\ell}(\tau)}{\mu(\tau)} \lambda(\tau)^{n}e^{-z\lambda(\tau)}.
\end{equation}
Then, we substitute Equation (\ref{eq:ftrace2}) in Equation (\ref{eq:dynamical-det}) and
we get that $\Cdet_{\ell}(\xi, z)$ can be interpreted as the ``determinant'' of
$(\COne - \xi R^{(\ell)}(z))^{-1}$, while $\Ddet_{\ell}(z)$ can be interpreted as 
the ``determinant'' of $z\COne - X^{(\ell)}$. Thus,
\begin{equation}
  \label{eq:ddet-ftrace}
  \Cdet(\xi - z, z) = \exp\left(-\sum_{i=0}^{\infty} \frac{(\xi - z)^{n}}{n}
  \ftrace(R^{(\ell)}(\xi)^{n})\right) 
  = \left( \prod_{\lambda_{i, \ell} \in B_{}(\xi, \rho_{p, q, \ell})}
  \frac{z - \lambda_{i, \ell}}{\xi - \lambda_{i, \ell}}  \right) \psi(\xi, z)
\end{equation}
where $\psi(\xi, z)$ is analytic and nonzero for $z \in B_{}(\xi, \rho_{p, q, \ell})$. \\ \\
Furthermore, Theorem \ref{theo:ddet-nonzero} implies Theorem \ref{theo:analytic-nonzero}. \qed

\begin{theo}[{\cite[Corollary 2.7]{GiLiPo}}]
  \label{theo:corre-aflow}
  The geodesic flow $\varphi_{t}: M \longrightarrow M$ for a compact manifold $M$
  with better than $\frac{1}{9}$\---pinched negative section curvatures is exponentially
  mixing with respect to the Bowen\---Margulis measure $\mu$; that is; there exists $\alpha$
  such that for $f, g \in C^{\infty}(T_{1}M)$ there exists a $C > 0$ for which the 
  correlation function
  $$
  \rho(t) = \int f \circ \varphi_{t}g d\mu - \int f d\mu \int g d\mu,
  $$
  satisfies $|\rho(t)| \leq C_{\#} e^{-\alpha |t|}$, for all $t \in \RR$.
\end{theo}
\proof Consider the Fourier transform $\fou(s) = \int_{-\infty}^{\infty} e^{ist} \rho(t) dt$
of the correlation function $\rho(t)$. By Theorem \ref{theo:po85} and \cite[Theorem 4.1]{Ru87},
the analytic extension of $\zruelle(z)$ in Theorem \ref{theo:cflow-pole}
implies that there exists $0 < \eta \leq \tau_{*}$ such that $\fou(s)$ has an analytic
extension to a strip $|\Imag(s)| < \eta$. Now, without the loss the generality,
we fixed each value $- \eta < t < \eta$, such we have that the function $\sigma \mapsto \fou
(\sigma + it)$ is in $L^{1}(\RR)$. Finally we apply the Paley\---Wiener Theorem \ref{theo:resi}
and result follows. \qed   \\ \\
As related to the Equation (\ref{eq:fourier}), the poles in $\fou(s)$ of Theorem 
\ref{theo:corre-aflow} are the Pollicott\---Ruelle resonances on a compact manifold
which asserts the assumptions $\hyperref[prop-B]{(B1)}$\---$\hyperref[prop-B]{(B2)}$.

\subsection{A short microlocal proof}
Unlike \cite{GiLiPo}; whose proofs only work for contact flows;
Dyatlov-Zworski \cite{DyZw13} proved in 2013, the meromorphic continuation of the Ruelle
zeta function for $C^{\infty}$ Anosov flows under the perspective of microlocal analysis, 
using semiclassical and scattering tools, and based on the study of the generator 
of the flow as a
semiclassical differential operator. The proofs applies to any Anosov flow for
which linearized Poincar\'e maps $\PC_{\gamma}$, where $\gamma$ is a closed orbit such that
$$
|\det(I - \PC_{\gamma})| = (-1)^{q}\det(I - \PC_{\gamma}), \hsp \text{ with } q 
\text{ independent of } \gamma.
$$
Furthermore, the assumptions $\hyperref[prop-B]{(B1)}$\---$\hyperref[prop-B]{(B2)}$ still hold
in this subsection. Let us first list some important definitions, for a major literature 
\--- see more
on microlocal analysis \cite{HoI-II, HoIII-IV, Ve, Mea, Iv}, semiclassical analysis 
\cite{Zwa, GuSt, EvZw, Be} and scattering theory \cite{Meb}.

%\begin{defi}[Dyatlov-Zworski \cite{DyZw15}] The Pollicott-Ruelle resonances are
%  limits (with multiplicities) of the eigenvalues $P + i\epsilon \Delta_{g}, \epsilon \rightarrow 0+$, where $-\Delta_{g} \leq 0$ is a Laplacian for some Riemannian metric $g$ on $X$.
%\end{defi}

\subsubsection{Definitions}
Victor Guillermin \cite{Gu}, defined a Trace formula using distributional operations
of pullback by some $\iota(t, x) = (t, x, x)$ and some pushforward $\pi:(t, x) \rightarrow t$
such that
$$
\ftrace e^{-itP} := \pi_{*}\iota^{*} K_{e^{-itP}},
$$
where $K_{\bullet}$ denotes the distributional kernel of operator \cite[Theorem 6]{Gu}.
As Lars H\"ormander claimed in \cite[Theorem 8.2.4]{HoI-II}, the pullback is well\---defined
in the sense of distributions since
  \begin{equation}
    \WF(K_{e^{-itP}}) \cap N^{*}(\RR_{t} \times \Delta(X)) = \emptyset, \hspace{1em} t > 0,
  \end{equation}
where $\Delta(X) \subset X \times X$ is the diagonal and $N^{*}(\RR_{t} \times \Delta(X))
\subset T^{*}(\RR_{t} \times X \times X)$ is the conormal bundle. Thus, we define
\begin{defi}[Guillemin's Trace Formula]
  \label{defi:gu-trace-formula}
  \begin{equation}
    \label{eq:gu-trace-formula}
    \ftrace e^{-itP} = \sum_{\gamma} \frac{T_{\gamma}^{\sharp} \delta(t - T_{\gamma})}
    {|\det(I - \PC_{\gamma})|}, \hspace{1em} \text{ for some }t > 0,
  \end{equation}
  where $T_{\gamma}$ is the period of the orbit $\gamma$, $T_{\gamma}^{\sharp}$ is the primitive
  period, $\PC_{\gamma}$ is the linearized Poincar\'e map and $\delta(\bullet)$
  is the Dirac delta function.
\end{defi}
 
For a detailed proof of Equation (\ref{eq:gu-trace-formula}) \--- see \cite[\ts Appendix B]{DyZw13}
and \cite[\ts II]{Gu}.  
Let $\WF(u)$ be the wavefront set for some $u \in \DC(M)$ distribution. Since we do need
a more robust measure of semiclassical regularity of functions, we define
the semiclassical wavefront set $\WF_{h}$ in the sense of \cite[\ts 8.4.2]{Zwa}, for a parameter
$h$ such $h$-tempered families of distributions $\{u(h)\}_{0 < h < 1}$.
\begin{defi}
  \label{defi:semi-wf}
  The semiclassical wavefront set $\WF_{h} \subset \overline{T}^{*}M$ is a subset from
  the fiber-radially compactified contangent\---bundle 
  (\ie a manifold with interior $T^{*}M$ and boundary $\partial \overline{T}^{*}M = 
  S^{*}M = (T^{*}M \backslash 0)/\RR^{+}$, the cosphere bundle \--- see \cite[\ts E.1]{DyZw}).
 Furthermore, $\WFh$ measures oscillations on the $h$\---scale and if $u$ is an $h$\---independent
  distribution, then 
  \begin{equation}
    \label{eq:-wf-relation}
    \WF(u) = \WF_{h}(u) \cap (T^{*}M \backslash 0).
  \end{equation}
\end{defi}  
Now, we consider the semiclassical operator $\PB \in \Psi_{h}^{k}(M, \Hom(\EC))$ where
$\EC$ is a vector bundle over $M$ such that it is acting on $h$\---tempered families
of distributions $\ub(h) \in \DC(M, \EC)$. From Definition \ref{defi:semi-wf}, we
denote the natural projection 
\begin{equation}
  \label{eq:nat-proj-semi}
  \kappa:T^{*}M \longrightarrow S^{*}M = \partial \overline{T}^{*}M.
\end{equation}  
Let $L \subset T^{*}M$ be a closed conic invariant set under the flow $e^{tHp}$
such there is an open neighbourhood $U$ of $L$ \--- see more \cite[\ts 18.3]{HoI-II}.
Then, the Equation (\ref{eq:nat-proj-semi}) asserts that
\begin{equation}
  \label{eq:prop-disc}
  \begin{split}
      d(\kappa(e^{-tH_{p}}(U)), \kappa(L)) \to 0 \quad &\text{ as } \quad t \to +\infty;  \\
      (x, \xi) \in U \implies |e^{-tH_{p}}(x, \xi)| \geq C^{-1} e^{\theta t}|\xi|,
      &\text{ for any norm on the fibers and some } \theta > 0.
  \end{split}
\end{equation}

\begin{defi}
  \label{defi:radial-source}
  Let $\kappa$ as in Equation (\ref{eq:nat-proj-semi}) and $L$ be a closed
  conic invariant that asserts Equation (\ref{eq:prop-disc}). Then, we say
  that $L$ is called a radial source and if we reverse the direction of 
  the flow, then $L$ is called a radial sink.
\end{defi}   
By Definition \ref{defi:radial-source} and letting $E_{s}^{*}, E_{0}^{*}, E_{u}^{*}$ be 
the duals of $E_{s}, E_{0}$ and $E_{u}$, respectively, then by Equation
(\ref{eq:prop-disc}), we say that $E_{s}^{*}$ and $E_{u}^{*}$ are a radial
source and a radial sink, respectively.  
Let $\PB$ as before such that $\PB: C^{\infty}(M; \EC) \longrightarrow 
C^{\infty}(M; \EC)$, besides 
$$
\PB(\ub) = \frac{1}{i} \Lop_{v}\ub, \hsp \EC = \bigoplus_{j=0}^{n} \Lambda^{j}(T^{*}M),
$$
where $V$ is the generator of the flow $\varphi_{t}$, $\Lop$ denotes
the Lie derivate and $\ub$ is a differential form on $M$.

\begin{defi}[Anisotropic Sobolev Spaces]
  \label{defi:ans-sob-space}
  The Anisotropic Sobolev spaces are defined using the exponential weight
  \--- see \cite[Lemma 7.6]{Zwa} and \cite[Theorem 7.7]{Zwa}
  \begin{equation}
    \label{eq:ans-sob-space}
    H_{sG} := \exp(-sG)(L^{2}(M)), \hspace{1em} \norm{\ub}_{H_{sG}} :=
    \norm{\exp(sG)\ub}_{L^{2}},
  \end{equation}
  where $G \in \Psi^{0+}(M)$ satisfying
  $$
  \sigma(G)(x, \xi) = m_{G}\log|\xi|,
  $$
  where $m_{G} = 1$ near $E_{s}^{*}$ and $m_{G} = -1$ near $E_{u}^{*}$.
\end{defi}  
For more about Anisotropic Sobolev spaces \--- see Duistermaat \cite{Du},
Unterberger \cite{Un}, Zworski \cite[\ts 8.3]{Zwa} and Baladi\---Tsujii \cite{BaTs07}.

\subsubsection{On the proof}
Now, let us use some essentials theorems from \cite{DyZw13} in order
to prove the existence of Pollicott\---Ruelle resonances on the compact case
via microlocal analysis.

\begin{theo}[\cite{DyZw13}]
  \label{theo:zruelle-theo}
  Supposse $M$ is a compact manifold and $\varphi_{t}:M \longrightarrow M$ is a $C^{\infty}$
  Anosov flow with orientable and unstable bundles. Let $\{\gsharp\}$ denote
  the set of primitive orbits of $\varphi_{t}$, with $T_{\gamma}^{\sharp}$ their periodics.
  Then the Ruelle zeta function,
  \begin{equation}
    \label{eq:rzeta-semi}
    \zruelle(\lambda) = \prod_{\gsharp}(1 - e^{i \lambda T_{\gamma}^{\sharp}}),
  \end{equation}
  which converges for $\Imag \lambda \gg 1$, has a meromorphic continuation to $\CC$.
\end{theo}  
One of the main Propositions in \cite{DyZw13}, is:

\begin{theo}[{\cite[Proposition 3.4]{DyZw13}}] 
  \label{theo:prop3.4}
    Fix a constant $C_0>0$ and $ \varepsilon > 0 $. Then for $s>0$ large
    enough depending on $C_0$
    and $h$ small enough, the operator
    $$
    \mathbf P_\delta ( z ) :D_{sG(h)}\to H_{sG(h)},\quad
    $$
    $$
    -C_0h\leq \Im z\leq 1,\quad
    |\Re z|\leq h^{\varepsilon},
    $$
    is invertible, and
    the inverse, $\mathbf R_\delta(z)$, satisfies
    $$
    \|\mathbf R_\delta(z)\|_{H_{sG(h)}\to H_{sG(h)}}\leq Ch^{-1},\quad
    $$
    $$
    \WFh'(\mathbf R_\delta(z))\cap T^*(M\times M)\subset \Delta(T^*M)\cup\Omega_+,
    $$
    with $\Delta(T^*M),\Omega_+$ defined in \cite[Propostion 3.3]{DyZw13},
    and %the semiclassical wavefront set 
    $\WFh'(
    %\mathbf R_\delta(z)
    \bullet)\subset\overline T^*(M\times M)$
    is defined for an $h$\---tempered family of operators $\bullet(h):C_{c}^{\infty}(M)
    \longrightarrow \DC(M)$.
\end{theo}
\proof The proof of Theorem \ref{theo:prop3.4} is assumed by $\norm{\ub}_{H_{sG(h)}}
\leq 1$ such that
\begin{equation}
  \label{eq:main-bound}
  \norm{\ub}_{H_{sG(h)}} \leq Ch^{-1}\norm{\fb}_{H_{sG(h)}}, \hsp 
  \ub \in D_{sG(h)}, \hsp \fb = \PB_{\delta}(z)\ub.
\end{equation}
Then for some $A \in \Psi_{h}^{0}(M)$, we can get bounds on $A\ub$ as are
detailed in \cite[Proposition 3.4]{DyZw13} which arrive to the Equation (\ref{eq:main-bound}).
\qed   \\ \\
From Theorem \ref{theo:prop3.4}, we can deduce that:
\begin{enumerate}[(1)]
\ii $H_{sG}$ and $D_{sG}$ are topologically
isomorphic to $H_{sG(h)}$ and $D_{sG(h)}$, respectively. And $Q_{\delta}:D_{sG} \longrightarrow H_{sG}$ is
smoothing and thus compact \--- see \cite[Proposition 3.1]{DyZw13}).
\ii If $\Imag \lambda > C_{1}, \ub \in H_{sG} \subset H^{-s}$ and $(\PB - \lambda)\ub =
\fb \in H_{sG}$, then
$$
\ub = - \int_{0}^{\infty} \partial_{t}(e^{i\lambda t} \varphi_{-t}^{*}\ub)dt
= i\int_{0}^{\infty} e^{i\lambda t} \varphi_{-t}^{*} \fb dt,
$$
where the integrals converge in $H^{-s}$. This also implies that $(\PB - \lambda)$ is
injective and invertible $D_{sG} \longrightarrow H_{sG}$. Then
\begin{equation}
  \label{eq:pr-resonances-micro}
    (\PB - \lambda)^{-1} = i \int_{0}^{\infty} e^{i\lambda t} \varphi_{-t}^{*} dt,
\end{equation}
where $\varphi_{-t}^{*}:C^{\infty}(X; \EC) \longrightarrow C^{\infty}(X; \EC)$ is
the pullback operator by $\varphi_{-t}$ on differential forms and the
integral on the right\---hand side converges in operator norm $H^{s} \rightarrow H^{s}$
and $H^{-s} \rightarrow H^{-s}$  \--- see \cite[Proposition 3.2]{DyZw13}.
\ii By \cite[\ts D.3]{Zwa}, $\RB(\lambda) = \RB_{H}(\lambda) + \sum_{j=1}^{J(\lambda_{0})}
A_{j}/(\lambda - \lambda_{0})^{j}$ where $\lambda_{0}$ is a near pole and $A_{j}$
are operators of finite rank such that
$$
\Pi := -A_{1} = \frac{1}{2\pi i} \oint_{\lambda_{0}} (\lambda - \PB)^{-1} d\lambda,
$$
where $[\Pi, \PB] = 0$. Thus $A_{j} = -(\PB - \lambda_{0})^{j-1} \Pi$ and
$(\PB - \lambda_{0})^{J(\lambda_{0})}\Pi = 0$ \--- 
see \cite[Proposition 3.3]{DyZw13}.
\ii Since $Q_{\delta}$ is pseudodifferential and supposing the fact that
$$
\RB(\lambda) = h(\RB_{\delta}(z) - i\RB_{\delta}(z)Q_{\delta}\RB_{\delta}(z))
- \RB_{\delta}(z)Q_{\delta}\RB(\lambda)Q_{\delta}\RB_{\delta}(z),
$$
we get that
$$
\WF_{h}'(\RB_{\delta}(z) - i\RB_{\delta}(z)Q_{\delta}\RB_{\delta}(z)) \cap
T^{*}(M \times M) \subset \Delta(T^{*}M) \cup \Omega_{+}
$$
\--- see \cite[Proposition 3.3]{DyZw13}.
\end{enumerate}
The Pollicott\---Ruelle resonances are the poles of $\Real(\lambda)$ in the
region $\Imag \lambda > -C_{0}$ of the meromorphic continuation of the
Schwartz Kernel of the operator given by the right\---hand side of (\ref{eq:pr-resonances-micro}),
and thus are independent of the choice of $s$ and the weight $G$.
 For the microlocal proof of the meromorphic continuation of Theorem \ref{theo:zruelle-theo}
\--- see \cite[\ts 4]{DyZw13}.

\subsubsection{Further developments}
Many applications had been development since the proof of those methods:

\begin{itemize}
    \ii Dyatlov\---Zworski \cite{DyZw15}, showed that Pollicott\---Ruelle resonances
    are the limits of eigenvalues of $V/i + i\varepsilon \delta_{g}$, as $\varepsilon
    \to 0+$, where $-\delta_{g}$ is any Laplace\---Beltrani operator on $X$.
    \ii Jin\---Zworski \cite{JiZw}, proved that for any Anosov flows there exists a strip with
    infinitely many resources and a counting function which cannot be sublinear.
    \ii Colin Guillarmou \cite{G1}, studied regularity properties of cohomological
    equations and provides applications. Guillarmou \cite{G2} also established a deformation
    lens rigidity for a class of manifolds including manifolds with negative curvature 
    and strictly convex boundary.
    \ii Dyatlov\---Guillarmou \cite{DyGu14}, proved meromorphic continuation for $(\PB - \lambda)^{-1}$
    and zeta functions for non\---compact manifolds with compact hyperbolic trapped sets.
    \ii Dyatlov\---Faure\---Guillarmou \cite{DyFaGu}, described the complex poles of the power
    spectrum of correlations for the geodesic flow on compact hyperbolic manifolds
    in terms of eigenvalues of the Laplacian on certain natural tensor bundles.
    \ii Dyatlov used \cite[Proposition 2.4]{DyZw13} and \cite[Proposition 2.5]{DyZw13}
    in \cite{Dya} as part to establish a resonance free strip for condimension 2
    symplectic normally hyperbolic trapped sets. To see a major literature about resonances
    for infinite\---area hyperbolic surfaces \--- see \cite{Bo16}.
    \ii Dyatlov\---Zworski \cite{DyZw15}, used microlocal methods similar to \cite{DyZw13}
    in order to show stochastic stability of Pollicott\---Ruelle resonances, more precisely,
    let $\PB_{\EC} = \frac{1}{i}V + i\EC\Delta_{g}$ and let $\{\lambda_{j}(\EC)\}_{0}^{\infty}$
    be the set of its $L^{2}$\---eigenvalues. Furthermore, let $\{\lambda_{j}\}_{j=0}^{\infty}$
    be the set of the Pollicott\---Ruelle resonances of the flow $\varphi_{t}$, then
    $\lambda_{j}(\EC) \to \lambda_{j}$ as $\EC \to 0+$ with convergence uniform
    for $\lambda_{j}$ in a compact set \--- see the proof in \cite[\ts 5]{DyZw15}.
    \ii Similar to \cite{DyZw15}, Zworski \cite{Zwc} showed scattering resonances
    of $-\Delta + V$ where $V \in L_{c}^{\infty}(\RR^{n})$, are the limits eigenvalues of
    $-\Delta + V - i\varepsilon x^{2}$ as $\varepsilon \to 0+$ via complex scaling
    method \cite[\ts 2]{Zwc} \--- to see more about scattering resonances \cite{DyZw, Zwb}.
    \ii Alexis Drouot \cite{Dr}, showed that for a compact manifold and negatively
    curved $\mathbb{M}$, the $L^{2}$\---spectrum of the infinitesimal generator of the
    Kinetic Brownian motion on the cosphere bundle as a stochastic process modeled
    by the geodesic equation perturbed with a random force of size $\varepsilon$,
    converges to the Pollicott\---Ruelle resonances as $\varepsilon$ goes to $0$.
\end{itemize}

\section{On the open systems case}
\label{section-open-case}
\hspace{-1.4em}
In 2014, Dyatlov\---Guillearmou \cite{DyGu14} defined Pollicott\---Ruelle resonances for
open systems, more precisely, geodesic flows on noncompact
asymptotically hyperbolic negatively curved manifolds, as well as for more general
open hyperbolic systems related to Axiom A flows. They used many generalized 
microlocal tools from \cite{DyZw13, FaSj} and functional analysis tools from \cite{GiLiPo},
and used anisotropic Sobolev spaces to control the singularities
at fiber infinity, and using complex absorbing potentials on the boundary and complex
absorbing pseudodifferential operators beyond the boundary to obtain a global Fredholm
problem for the extension of $\XB$ to a compact manifold without boundary $\MC$ \--- see more
\cite[\ts 4]{DyGu14}.

\subsection{Definitions}
We use the same notation as in \cite{DyGu14}, given a $n$\---dimensional compact manifold
$\overline{\UC}$ with interior $\UC$ and boundary $\partial \UC$, then $X$ is a smooth
$C^{\infty}$\---nonvanishing vector field on $\overline{\UC}$ such that for some $t$, the
corresponding flow is defined as $\varphi^{t} = e^{tX}$. Furthermore, $\partial\UC$
is strictly convex (\ie for some $x \in \UC$ then $X \rho(x) = 0 \implies X^{2}\rho(x) < 0$
where $\rho \in C^{\infty}(\overline{\UC})$).

\begin{defi}
  \label{defi:in-out-tails}
  The incoming ($\Gamma_{+}$) and outgoing ($\Gamma_{-}$) tails are subsets from $\overline{\UC}$
  such that
  \begin{equation}
    \label{eq:in-out-tails}
    \Gamma_{\pm} = \bigcap_{\pm t \geq 0} \varphi^{t}(\overline{U}).
  \end{equation}
\end{defi}  
From Definition \ref{defi:in-out-tails}, let $K = \gplus \cap \gminus$ be the trapped set
such that for some $x \in K$ there is a splitting in $T_{x}\MC$ in the sense of Equation
(\ref{eq:splitting}) and Equation (\ref{eq:hyper}), where $\MC$ is a compact
manifold without boundary such that $\overline{\UC}$ is embedded in $\MC$ \--- see
more about dynamical assumptions of the trapped set in \cite[\ts 3.5.1]{Dyb}. 
Now, let $\EC$ be the smooth complex vector bundle over $\overline{\UC}$ and the
first order differential operator $\XB: C^{\infty}(\overline{\UC}; \EC) \longrightarrow
C^{\infty}(\overline{\UC}; \EC)$ such that
\begin{equation}
  \label{eq:diff-op}
  \XB(f\ub) = (Xf)\ub + f(\XB\ub), \hsp f \in C^{\infty}(\overline{U}), 
  \ub \in C^{\infty}(\overline{\UC}; \EC).
\end{equation}
Fixing a smooth measure $\mu$ on $\MC$ and the norm $L^{2}(\MC; \EC)$, we
define the transfer operator $e^{-t\XB}:L^{2}(\MC; \EC) \longrightarrow L^{2}(\MC; \EC)$
and by Equation (\ref{eq:diff-op}) we have that the support of $e^{-t\XB}$ is:
\begin{equation}
  \label{eq:support-e}
  e^{-t\XB}(f\ub) = (f \circ \varphi^{-t})e^{-t\XB}\ub, \hsp f \in C^{\infty}(\MC),
  \ub \in C^{\infty}(\MC; \EC).
\end{equation}

\begin{defi}
  \label{defi:resonant-state}
  Let $\RB$ be the restricted resolvent defined as
  \begin{equation}
    \label{eq:family-res}
      \RB(\lambda) = \COne_{\UC}(\XB + \lambda)^{-1}\COne_{\UC}:C^{\infty}_{0}(\UC; \EC)
      \longrightarrow \DC(\UC; \EC), \hsp \Real \lambda > C_{0},
  \end{equation}
  where $\lambda \in \CC$. Then, for each $j \geq 1$ the space
  of generalized resonant states
  \begin{equation}
    \label{eq:gen-resonant-state}
    \Res_{\XB}^{(j)}(\lambda) = \{\ub \in \DC(\UC; \EC) \ | \ 
    \supp \ub \subset \gplus, \WF(\ub) \subset E^{*}_{+}, 
  (\XB + \lambda)^{j}(\ub) = 0 \},
  \end{equation}
  where $E^{*}_{+} \supset E^{*}_{u}$ is the extended unstable bundle over
  $\gplus$.
\end{defi} \hspace{-1.4em}
For a detailed construction of $E_{+}^{*}$ in Definition \ref{defi:resonant-state},
see \cite[Lemma 2.10]{DyGu14}. Furthermore, the subbundle $E_{+}^{*}$ is a
generalized radial sink and $E_{-}^{*} \supset E_{s}^{*}$ is a generalized
radial source; which is a modification from Equation (\ref{eq:prop-disc}).  
As related to Definition \ref{defi:ans-sob-space}, \cite{DyGu14} defined the anisotropic 
Sobolev space $\HC_{h}^{r}$; in order to control the singularities at fiber infinity; as
\begin{equation}
  \label{eq:ans-sob-h-open}
  \HC_{h}^{r} = \exp(-rG(h))(L^{2}(\MC; \EC)), \hsp \norm{\ub}_{\HC_{h}^{r}} =
  \norm{\exp(rG(h))\ub}_{L^{2}(\MC; \EC)},
\end{equation}
where $G$ is the operator defined as $G(h) \in \bigcap_{\lambda > 0} \Psi_{h}^{k}(\MC)$
\--- see more \cite[\ts 4.1]{DyGu14} and for the propagation of singularities
\cite[Proposition 2.5]{DyZw13}.  
Now, let $V \in C^{\infty}(\UC; \CC)$, then $\gamma^{\sharp}:[0, T_{\gamma^{\sharp}}]
\longrightarrow K$ of $\varphi_{t}$ of period $T_{\gamma^{\sharp}}$, thus
\begin{equation}
  \label{eq:gsharp-open}
  V_{\gamma^{\sharp}} = \frac{1}{T_{\gamma^{\sharp}}} 
  \int_{0}^{T_{\gamma^{\sharp}}} V(\gamma^{\sharp}(t)) dt
\end{equation}
be the average of $V$ over $\gamma^{\sharp}$. Thus, we define the Ruelle zeta function
as the product over all primitive closed trajectories of $\varphi^{t}$ on $K$:

\begin{equation}
  \label{eq:zruelle-open}
  \zeta_{\text{Ruelle} V}(\lambda) = \prod_{\gamma^{\sharp}}(1 - \exp(-T_{\gamma}^{\sharp}(
  \lambda + V_{\gamma^{\sharp}}))), \hsp \Real \lambda \gg 1.
\end{equation}
For express Pollicott\---Ruelle resonances of $\XB$ as poles, let $\EC_{0}$ be
the vector bundle over $\overline{\UC}$ by
$$
\EC_{0}(x) = \{\eta \in T_{x}^{*}\MC \ | \ \iprod{X(x)}{\eta} = 0\}, \hsp
x \in \overline{\UC},
$$
and let $\PC_{x, t}:\EC_{0}(x) \longrightarrow \EC_{0}(\varphi^{t}(x))$ be
the Poincar\'e map such that $\PC_{x, t} = (d\varphi^{t}(x))^{-T}|_{\EC_{0}(x)}$.
Now, for each $\ub \in C^{\infty}(\MC; \EC)$, we put $\alpha_{x, t}(\ub(x)) =
e^{-t\XB}\ub(\varphi^{t}(x))$ where $\alpha_{x, t}$ is the parallel transport
defined as $\alpha_{x, t}: \EC(x) \longrightarrow \EC(\varphi^{t}(x))$, then if
$\ub(x) = 0$ implies that $e^{-t\XB}\ub(\varphi^{t}(x)) = 0$ by Equation (\ref{eq:support-e}).
Thus, for the operator $\alpha_{\varphi^{t}(x_{0}),  T}:\EC(\varphi^{t}(x_{0})) 
\longrightarrow \EC(\varphi^{t}(x_{0}))$ where $T > 0$ we have that:
\begin{equation}
  \label{eq:tr-open}
  \tr \alpha_{\varphi^{t}(x_{0})} = \tr \alpha_{\varphi^{t}(x_{0}), T}, \hsp
  \det(I - \PC_{\gamma}) = \det(I - \PC_{\varphi^{t}(x_{0}), T}) \neq 0.
\end{equation}

\begin{defi}
Using the wavefront set $\WF$ in the sense of Definition \ref{defi:semi-wf},
of any $u \in \DC(\MC)$ and considering wavefront sets $\WF'(B) \subset T^{*}
(\MC \times \MC) \backslash 0$, where $B:C^{\infty}(\MC) \longrightarrow \DC(\MC)$
are operators, we define 
\begin{equation}
  \label{eq:wf2}
  \WF'(B) = \{(x, \xi, y, -\eta) \ | \ (x, \xi, y, \eta) \in \WF(K_{B})\},
\end{equation}
where the Schwartz Kernel $K_{B} \in \DC(\MC \times \MC)$ is given by
\begin{equation}
  Bf(x) = \int_{\MC} K_{B}(x, y)f(y)dy, \hsp f \in C^{\infty}(\MC).
\end{equation}
\end{defi}  
Let $V, W \in \overline{T}^{*}\MC$ be open sets, such that $e^{-TH_{p}}(x, \xi) \in V$
and $e^{-tH_{p}}(x, \xi) \in W$ for $t \in [0, T]$. We denote the open subset
$$
\Con_{p}(V; W) \subset \overline{T}^{*}\MC,
$$
the set of such points \--- see \cite[Proposition 2.5]{DyGu14}. Let
$A, B, B_{1} \in \Psi_{h}^{0}(\MC)$ be operators such that $q \geq 0$ near $\WF_{h}(B_{1})$,
where $\WF_{h} \subset \Con_{p}(\Ellh(B); \Ellh(B_{1}))$, thus the trajectories
of $e^{-tH_{p}}$ starting on $\WF_{h}(A)$ either pass though $\Ellh(B)$ or
converge to some closed set $L$, while staying on $\Ellh(B_{1})$ 
\--- see \cite[Definition 3.3]{DyGu14}. 

\subsection{On the proof}
Dyatlov\---Guillarmour used sharp G\r{a}rding inequatlity \--- see \cite[\ts 4.7]{Zwa},
in \cite[Lemma 3.4]{DyGu14} to show that the definition of real part $\Re \PB$ is not
trivial, that is, assume that $\PB \in \Psi_{h}^{2m + 1}(\MC; \EC)$ is principally
scalar, $A \in \Psi_{h}^{0}(\MC)$, and $\Re \sigma_{h}(\PB) \leq 0$ in a neighborhood
$U \subset \CPC$ of $\WF_{h}(A)$. Then, there exist a constant $C$ such that for
each $N$ and $\ub \in H_{h}^{m+1/2}(\MC, \EC)$,
\begin{equation}
  \label{eq:per1}
  \Re\iprod{\PB A\ub}{A\ub}_{L^{2}} \leq Ch\norm{A\ub}_{H_{h}^{m}}^{2} + 
  \mathcal{O}(h^{\infty})\norm{\ub}_{H_{h}^{-N}}^{2}.
\end{equation}
Furthermore, if $L$ and $\PB$ satisfies that $\Im(\mathbf P-iQ)\lesssim -h$
on $H^s_h$ near $L$, for all $s$, where $Q \in \Psi_{h}^{0}(\MC)$,
$\Im\sigma_h(\mathbf P)\leq 0$ near $L$ and $\Re\sigma_h(Q)>0$ on $L$. Then, for some
aditional $p:=\Re\sigma_h(\mathbf P)\in\Hom^1(T^*\mathcal M;\mathbb R)$
and assuming that $L\subset \partial\overline T^*\mathcal M$, where
$L$ is invariant under~$e^{tH_p}$. Fix a metric $|\cdot|$ on the fibers of $T^*\mathcal M$.
Then,
\begin{enumerate}
    \ii Assume that there exist $c, \gamma > 0$ such that
    \begin{equation}
    \label{eq:per2}
      \frac{|e^{tH_{p}(x, \xi)}|}{|\xi|} \geq ce^{\gamma|t|} \quad \text{for } 
      (x, \xi) \in L, t \leq 0. 
    \end{equation}
    Then there exists $s_{0}$ such that for all $s > s_{0}$, $\Im \PB \lesssim -h$
    near $L$ on $H_{h}^{s}$.
    \ii Assume that there exist $c, \gamma > 0$ such that
    \begin{equation}
    \label{eq:per3}
      \frac{|e^{tH_{p}(x, \xi)}|}{|\xi|} \geq ce^{\gamma|t|} \quad \text{for } 
      (x, \xi) \in L, t \geq 0. 
    \end{equation}
    Then there exists $s_{0}$ such that for all $s < s_{0}$, $\Im \PB \lesssim -h$
    near $L$ on $H_{h}^{s}$.
\end{enumerate}  
For the proofs of Equations (\ref{eq:per1}), (\ref{eq:per2}) and (\ref{eq:per3}) \--- see
\cite[\ts 3]{DyGu14}.

%\begin{theo}[{\cite[Theorem 3]{DyGu14}}]
%  Assume that the stable/unstable foliations $E_{u}, E_{s}$ are orientable. Then
%  the function $\zruelle_{V}(\lambda)$ admits a meromorphic continuation to 
%  $\lambda \in \CC$.
%\end{theo}
%\proof 

\begin{theo}
  \label{theo:cp-entire}
  The family of $\{\RB(\lambda)\}$, defined in the sense of Equation (\ref{eq:family-res}),
  continues meromorphically to $\lambda \in \CC$, with poles of finite rank.
\end{theo}

\begin{theo}[{\cite[Theorem 4]{DyGu14}}]
  \label{theo:main-theo}
  Define for $\Real \lambda \gg 1$
  \begin{equation}
    \label{eq:main-eq-open}
    F_{\XB}(\lambda) = \sum_{\gamma} \frac{e^{-\lambda T_{\gamma}} 
    T_{\gamma}^{\sharp}\tr \alpha_{\gamma}}{|\det(I - \PC_{\gamma})|},
  \end{equation}
  where the sum is over all closed trajectories $\gamma$ inside $K, T_{\gamma} > 0$
  is the period of $\gamma$, and $T_{\gamma}^{\sharp}$ is the primitive period.
  Then $F(\lambda)$ extends meromorphically to $\lambda \in \CC$. The poles
  of $F(\lambda)$ are the Pollicott\---Ruelle resonances of $\XB$
  and the residue at a pole $\lambda_{0}$ is equal to the rank
  of $\Pi_{\lambda_{0}}$.
\end{theo}
\proof We define the flat trace in the sense of the operator $\AB:C^{\infty}
(\MC; \UC) \longrightarrow \DC(\MC; \UC)$ such that $\WF'(\AB) \cap \Delta
(T^{*}\MC \backslash 0) = \emptyset$, then
\begin{equation}
  \label{eq:flat-open}
  \ftrace \AB = \int_{\MC} \tr_{\End(\EC)} K_{\AB}(x, x)dx.
\end{equation}
Making 
\begin{equation}
  \label{eq:fx}
F_{\XB}(\lambda) = \ftrace(\chi e^{-t_{0}(\XB + \lambda)}\RB(\lambda)\chi)
\end{equation}
for $\lambda > C_{1}$, $C_{1} > 0$, $\chi \in C_{0}^{\infty}(\UC)$ 
and $t_{0} > 0$ is small enough so that $t_{0} < T_{\gamma}$ for all $\gamma$.
Then, by \cite[Theorem 2]{DyGu14} and Equation (\ref{eq:fx}) we have that
\begin{equation}
  \label{eq:sum-trace-open}
  \tr^\flat\sum_{j=1}^{J(\lambda_0)}(-1)^{j-1}{\chi e^{-t_0(\mathbf X+\lambda)}
  (\mathbf X+\lambda_0)^{j-1}\Pi_{\lambda_0}\chi\over (\lambda-\lambda_0)^j}
  ={\rank\Pi_{\lambda_0}\over \lambda-\lambda_0}+\Hol(\lambda),
\end{equation}
where $\Hol(\lambda)$ is holomorphic near $\lambda_{0}$. \qed  \\ \\
The proof (and in fact the work of Dyatlov\---Guillarmou) is really complex, for
a full detailed and several particular cases of Theorem \ref{theo:main-theo} \--- see
\cite[\ts 4]{DyGu14}.  
By \cite[Lemma 4.3]{DyGu14} and \cite[Lemma 3.3]{DyGu14}, the operator
$-ih\COne_{\UC}\RB_{0}(ih\lambda)\COne_{\UC}$ gives the meromorphic continuation of $\RB(\lambda)$
in the region $[-C_{1}, h^{-1}] + i[-C_{2}, C_{2}]$ for $h$ small enough. Since
$C_{1}$ and $C_{2}$ can be chosen arbitraly and $h$ can be arbitrally small, we obtain the
continuation to the entire complex plane and Theorem \ref{theo:cp-entire} follows. \qed

%%%%%%%%%%%%%%%%%%%%%%%%%%%%%%%%%%%%%%%%%%%%%%%%%%%%%%%%%%%%%%%%%%%%%%%%%%%%%%%%
%%%%%%%%%%%%%%%%%%%%%%%%%%%%%%%%%%%%%%%%%%%%%%%%%%%%%%%%%%%%%%%%%%%%%%%%%%%%%%%%


\begin{thebibliography}{0}

\bibitem[AFW]{AFW} Jean-François Arnoldi, Fr\'ed\'eric Faure and Tobias Weich,
  \emph{Asymptotic spectral gap and Weyl law for Ruelle resonances of open partially
  expanding maps, \/}
  Erg. Thoery Dyn. Syst. \textbf{} (2015), 1-58; \arXiv{1302.3087}. 

\bibitem[Ba16]{Ba16} Viviane Baladi,
  \emph{Dynamical zeta Functions and Dynamical Determinants for Hyperbolic Maps
  \--- A Functional Approach, \/}
  book in preparation: \url{https://webusers.imj-prg.fr/~viviane.baladi/baladi-zeta2016.pdf}.

\bibitem[BaLi]{BaLi} Viviane Baladi and Carlangelo Liverani,
  \emph{Exponential decay of correlations for piecewise cone hyperbolic contact flows, \/}
  Comm. Math. Phys. \textbf{314} (2012), 689-773.
  
\bibitem[BaTs07]{BaTs07} Viviane Baladi and Masato Tsujii,
  \emph{Anisotropic H\"older and Sobolev spaces for hyperbolic diffeomorphisms, \/}
  Ann. Inst. Fourier \textbf{57} (2007), 127-154.

\bibitem[BaTs08]{BaTs08} Viviane Baladi and Masato Tsujii,
  \emph{Dynamical determinants and spectrum for hyperbolic diffeomorphisms, \/}
  Geometric and probabilistic structures in dynamics. Amer. Math. Soc. (2008), 29-68.

\bibitem[Be]{Be} Ivan Belisario Ventura,
  \emph{Applications of semiclassical analysis to partial differential equations, \/}
  UC Berkeley: Mathematics (2012). Retrieved from: \url{http://escholarship.org/uc/item/8jc5396t}.

\bibitem[Bo]{Bo} Rufus Bowen,
  \emph{Symbolic dynamics for hyperbolic flows, \/}
  Amer. J. Math. \textbf{95} (1973), 429-460.

\bibitem[Bo16]{Bo16} David Borthwick,
  \emph{Spectral theory of infinite-area hyperbolic surfaces, Second Edition, \/}
  Springer International Publishing, 2016.

\bibitem[BuLi]{BuLi} Oliver Butterley and Carlangelo Liverani,
  \emph{Smooth Anosov flows: Correlation spectra and stability, \/}
  J. Mod. Dyn. (1) \textbf{2} (2007), 301-322.

\bibitem[DDZ]{DDZ} Kiril Datchev, Semyon Dyatlov, Maciej Zworski,
  \emph{Sharp polynomial bounds on the number of Pollicott-Ruelle resonances, \/}
  Erg. Theory Dyn. Syst. \textbf{34} (2014), 1168-1183.

\bibitem[Dr]{Dr} Alexis Drouot,
  \emph{Pollicot\---Ruelle resonances via Kinetic Brownian motion}
  preprint, \arXiv{1607.03841}

\bibitem[Dya]{Dya} Semyon Dyatlov,
  \emph{Spectral gaps for normally hyperbolic trapping, \/}
  preprint \arXiv{1403.6401}.

\bibitem[Dyb]{Dyb} Semyon Dyatlov,
  \emph{Resonances in general relativity, \/}
  UC Berkeley: Mathematics. Retrieved from: \url{http://escholarship.org/uc/item/59d425g7}.

\bibitem[DyFaGu]{DyFaGu} Semyon Dyatlov, Fr\'ed\'eric Faure and Colin Guillarmou,
  \emph{Power spectrum of the geodesic flow on hyperbolic manifolds, \/}
  Analysis \& PDE \textbf{8} (2015), 923-2000.

\bibitem[DyGu14]{DyGu14} Semyon Dyatlov and Colin Guillarmou,
  \emph{Pollicott-Ruelle resonances for Open Systems, \/}
  C. Ann. Henri Poincar\'e, \textbf{11} (2016), 3089-3146; \arXiv{1410.5516}.

\bibitem[DyZw]{DyZw} Semyon Dyatlov and Maciej Zworski,
  \emph{Mathematical theory of scattering resonances, \/}
  book in preparation: \url{http://math.mit.edu/~dyatlov/res/}.

\bibitem[DyZw13]{DyZw13} Semyon Dyatlov and Maciej Zworski,
  \emph{Dynamical zeta functions for Anosov flows via microlocal analysis, \/}
  Annales de l'ENS \textbf{49} (2016), 543-577; \arXiv{1306.4203}.

\bibitem[DyZw15]{DyZw15} Semyon Dyatlov and Maciej Zworski,
  \emph{Stochastic stability of Pollicott-Ruelle resonances, \/}
  Nonlinearity, \textbf{28} (2015), 3511-3534;
  \arXiv{1407.8531}.

\bibitem[Du]{Du} Johannes Duistermaat,
  \emph{On Carleman estimates for pseudo-differential operatos, \/}
  Invent. Math. \textbf{17} (1972), 31-43.

\bibitem[EvZw]{EvZw} Lawrence Evans and Maciej Zworski,
  \emph{Lectures on Semiclassical Analysis, \/}
  notes in preparation: \url{https://math.berkeley.edu/~evans/semiclassical.pdf}

\bibitem[FaSj]{FaSj} Fr\'ed\'eric Faure and Johannes Sj\"ostrand,
  \emph{Upper bound on the density of Ruelle resonances for Anosov flow, \/}
  Comm. Math. Phys. \textbf{308} (2011), no. 2, 325-364; \arXiv{1003.0513}.

\bibitem[FaTsa]{FaTsa} Fr\'ederic Faure and Masato Tsujii,
  \emph{Band structure of the Ruelle spectrum of contact Anosov flows, \/}
  Comptes rendus \--- Math ́ematique \textbf{351} (2013), 385-391.

\bibitem[FaTsb]{FaTsb} Fr\'ederic Faure and Masato Tsujii,
  \emph{The semiclassical zeta function for geodesic flows on negatively curved manifolds, \/}
  To appear in Inventiones; \arXiv{1311.4932}.

\bibitem[Fr]{Fr} David Fried,
  \emph{Meromorphic zeta functions for analytic flows, \/}
  Comm. Math. Phys. (1) \textbf{174} (1995), 161-190.

\bibitem[G1]{G1} Colin Guillarmou,
  \emph{Invariant distributions and X-ray transform for Anosov flows, \/}
  preprint, \arXiv{1408.4732}.

\bibitem[G2]{G2} Colin Guillarmou,
  \emph{Lens rigidity for manifolds with hyperbolic trapped set, \/}
  preprint, \arXiv{1412.1760}.

\bibitem[GiLiPo]{GiLiPo} Paolo Giulietti, Carlangelo Liverani and Mark Pollicott,
  \emph{Anosov flows and dynamical zeta functions, \/}
  Ann. of Math. (2) \textbf{178} (2013), 687-773.

\bibitem[GoLi]{GoLi} S\'ebastien Gou\"ezel and Carlangelo Liverani,
  \emph{Banach spaces adapted to Anosov systems, \/}
  Erg. Theory Dyn. Syst. \textbf{26} (2006), 189-217.

\bibitem[Gu]{Gu} Victor Guillermin,
  \emph{Lectures on spectral theory of elliptic operators, \/}
  Duke Math. J. \textbf{44}(1977), 485-517.

\bibitem[GuSt]{GuSt} Victor Guillermin and Shlomo Stemberg,
  \emph{Semi-Classical Analysis, \/}
  Int. Press of Boston, 2013.

\bibitem[Ha]{Ha} Nicolai T. A. Haydn,
  \emph{Meromorphic extension of the zeta function for Axiom A flows, \/}
  Erg. Theory Dyn. Syst. \textbf{10} (1990), 347-360.

\bibitem[H\"oI-II]{HoI-II} Lars H\"ormander,
  \emph{The Analysis of Linear Partial Differential Operators, Volumes I and II, \/}
  Springer, 1983.

\bibitem[H\"oIII-IV]{HoIII-IV} Lars H\"ormander,
  \emph{The Analysis of Linear Partial Differential Operators, Volumes III and IV, \/}
  Springer, 1985.

\bibitem[Iv]{Iv} Victor Ivrii,
  \emph{Microlocal Analysis, Sharp Spectral Asymptotics and Applications, \/}
  book in preparation: \url{http://weyl.math.toronto.edu/victor_ivrii/research/future-book/}.

\bibitem[JiZw]{JiZw} Long Jin and Maciej Zworski,
  \emph{A local trace formula for Anosov flows,\/}
  with an appendix by Fr\'ed\'eric Naud,
  preprint, \arXiv{1411.6177}.

\bibitem[KaHa]{KaHa} Anatole Katok and Boris Hasselblatt,
  \emph{Introduction to the modern theory of dynamical systems, \/}
  Cambridge Univ. Press 1997.

\bibitem[Lia]{Lia} Carlangelo Liverani,
  \emph{On contact Anosov flows, \/}
  Ann. of Math. (3) \textbf{159} (2004), 1275-1312.

\bibitem[Lib]{Lib} Carlangelo Liverani,
  \emph{Decay of correlations, \/}
  Ann. of Math. \textbf{142} (1995), 239-301.

\bibitem[LiTs]{LiTs} Carlangelo Liverani and Masato Tsujii,
  \emph{Zeta functions and dynamical systems, \/}
  Nonlinearity, (19) \textbf{10} (2006), 2467-2473.

\bibitem[PaPo]{PaPo} William Parry and Mark Pollicott,
  \emph{Zeta functions and the periodic orbit structure of hyperbolic manifolds, \/}
  In \textit{Ast\'erisque} vol. 187-188. Soci\'et\'e math\'ematique de France, 1990.

\bibitem[Po85]{Po85} Mark Pollicott,
  \emph{On the rate of mixing of Axiom A flows, \/}
  Invent. Math. \textbf{81} (1985), 413-426.

\bibitem[Po86]{Po86} Mark Pollicott, 
  \emph{Meromorphic extensions of generalized zeta functions, \/}
  Invent. Math \textbf{85} (1986), 147-164.

\bibitem[ReSi]{ReSi} Michael Reed and Barry Simon,
  \emph{Functional Analysis, II, \/}
  Addison-Wesley, 1980.

\bibitem[Ru76]{Ru76} David Ruelle,
  \emph{Zeta-Functions for expanding maps and Anosov flows, \/}
  Invent. Math \textbf{34} (1976), 231-242.

\bibitem[Ru86]{Ru86} David Ruelle,
  \emph{Resonances of chaotic dynamical systems, \/}
  Phys. Rev. Lett. \textbf{56} (1986), 405-407.

\bibitem[Ru87]{Ru87} David Ruelle,
  \emph{Resonances for Axiom A flows, \/}
  J. Diff. Geom. \textbf{25} (1987), 99-116.

\bibitem[Mea]{Mea} Richard B. Melrose,
  \emph{Introduction to Microlocal Analysis, \/}
  book in preparation: \url{http://www-math.mit.edu/~rbm/iml90.pdf}

\bibitem[Meb]{Meb} Richard B. Melrose,
  \emph{Geometric Scattering Theory, \/}
  Cambridge University Press 1995.

\bibitem[Sj]{Sj} Johannes Sj\"ostrand,
  \emph{Geometric bounds on the density of resonances for semiclassical problems, \/}
  Duke Math. J. (1) \textbf{60} (1990), 1-57.

\bibitem[Sm]{Sm} Steven Smale, 
  \emph{Differentiable dynamical systems, \/}
  Bull. Amer. Math. Soc. \textbf{73} (1967), 747-817.

\bibitem[Un]{Un} Andr\'e Unterberger,
  \emph{R\'esolution d'\'equations aux d\'eriv\'ees partielles dans des espaces de
  distributions d'ordre de r\'egularit\'e, \/}
  Ann. Inst. Fourier \textbf{21} (1971), 85-128.

\bibitem[Na]{Na} Fr\'ed\'eric Naud,
  \emph{Analyse spectrale d'op\'erateurs de transfert et R\'esonances, \/}
  M\'emoire d'HDR, Math\'ematiques section 25;
  \url{http://fredericnaud.perso.sfr.fr/Habilitation.pdf}.

\bibitem[Ve]{Ve} Yves Colin de Verdi\`ere,
  \emph{M\'ethodes Semi-Classiques et Th\'eorie Spectrale, \/}
  book in preparation: \url{https://www-fourier.ujf-grenoble.fr/~ycolver/All-Articles/93b.pdf}

\bibitem[Wi]{Wi} Amie Wilkinson,
	\emph{Lectures on marked length spectrum rigidity, \/}
	lecture notes, \url{http://www.math.uchicago.edu/~wilkinso/papers/PCMI-Wilkinson.pdf}.

\bibitem[Zwa]{Zwa} Maciej Zworski,
	\emph{Semiclassical analysis, \/}
	Graduate Studies in Mathematics \textbf{138}, AMS, 2012.

\bibitem[Zwb]{Zwb} Maciej Zworski,
  \emph{Mathematical study of scattering resonances, \/}
  preprint, \arXiv{1609.03550}.

\bibitem[Zwc]{Zwc} Maciej Zworski,
	\emph{Scattering resonances as viscosity limits, \/}
  preprint, \arXiv{1505.00721}.

\end{thebibliography}
\end{document}